\documentclass[12pt]{article}

\usepackage{amsmath}
\usepackage{amsfonts}
\usepackage{amsthm}
\usepackage[colorinlistoftodos]{todonotes}
\usepackage{pifont}
\usepackage{mdframed,color}
\usepackage{hyperref}
\usepackage{pgfplots}
\pgfplotsset{compat=1.13}

\newenvironment{atmProof}{\noindent\ignorespaces\paragraph{Proof:}}{\hfill \ding{122}\par\noindent}
\newtheorem{theorem}{Theorem}
\pgfplotsset{
    midpoint segments/.code={\pgfmathsetmacro\midpointsegments{#1}},
    midpoint segments=3,
    midpoint/.style args={#1:#2}{
        ybar interval,
        domain=#1+((#2-#1)/\midpointsegments)/2:#2+((#2-#1)/\midpointsegments)/2,
        samples=\midpointsegments+1,
        x filter/.code=\pgfmathparse{\pgfmathresult-((#2-#1)/\midpointsegments)/2}
    }
}

\pgfplotsset{
    right segments/.code={\pgfmathsetmacro\rightsegments{#1}},
    right segments=3,
    right/.style args={#1:#2}{
        ybar interval,
        domain=#1+((#2-#1)/\rightsegments):#2+((#2-#1)/\rightsegments),
        samples=\rightsegments+1,
        x filter/.code=\pgfmathparse{\pgfmathresult-((#2-#1)/\rightsegments)}
    }
}

\pgfplotsset{
    left segments/.code={\pgfmathsetmacro\leftsegments{#1}},
    left segments=3,
    left/.style args={#1:#2}{
        ybar interval,
        domain=#1:#2,
        samples=\leftsegments+1,
        x filter/.code=\pgfmathparse{\pgfmathresult}
       }
}

\title{Infinitely Many Primes Using Generating Functions} 
\author{Sandeep Silwal \thanks{Email: \url{silwal@mit.edu}}}     
\begin{document}
\date{}
\maketitle

\section{Introduction}
Euclid gave the first proof that there are infinitely many prime numbers more than $2000$ years ago. His proof relied on the observation that if $p_1, \cdots, p_n$ are all prime numbers, then $p_i$ does not divide $Q := \left(\prod_{i=1}^n p_i \right) + 1$ and thus we can find a new prime among the prime factors of $Q.$ Therefore, $\mathbb P$, the set of prime numbers, must be infinite. Since then, many proofs have been given for this fact using techniques from a wide range of mathematics such as analysis, information theory, and even topology. We refer the reader to \cite{thebook, rosen} for a collection of such proofs.

Aside from Euclid's proof, perhaps the most well-known approach is Euler's proof that, because the series $\sum_{p \in \mathbb P} \frac{1}{p}$ diverges to infinity, the set of primes must be infinite (or else the sum would just be a finite rational number). Euler did not rigorously justify this proof, but future authors polished his argument. For a discussion on Euler's proof, see \cite{euler}.

Here we offer a proof of the infinitude of primes in the spirit of Euler's approach; we prove the divergence of the series 
\begin{equation} \label{first}
    \sum_{p \in \mathbb P} \frac{1}{\log(p)}, 
\end{equation}
which implies $\mathbb P$ is infinite in much the same way. In fact, because $\frac{1}{p} < \frac{1}{\log(p)}$ for all primes $p \geq 3$, the divergence of the sum of the reciprocals of $\log(p)$ follows immediately from Euler's proof. Instead on relying on Euler's proof, we use generating functions and tools from calculus to show that \eqref{first} diverges which implies that $\mathbb P$ is infinite. Specifically, we will prove the following theorem.

\begin{theorem}[Main Theorem] \label{thm} The following inequality holds for sufficiently large $n$:  $$\sum_{ p \le n } \frac{1}{\log(p)} > \frac{1}3 \log(n).$$ 
\end{theorem}
A corollary of Theorem \ref{thm} will be that there are infinitely many primes. We note that the estimate $\sum_{ p \le n} \frac{1}{\log(p)} = \Omega(\log(n))$ is very weak. In particular, using the prime number theorem, which states that the number of primes less than $x$ is asymptotically $\frac{x}{\log (x)}$ \cite{rosen, PNT}, we can actually prove the stronger estimate 
$$\sum_{p \le n} \frac{1}{\log(p)} = \Omega\left( \frac{n}{\log^2(n)} \right)$$
which we will not show here. Even though our approach yields much weaker results, our proof highlights the usefulness of generating functions and only uses elementary calculus. We now proceed to prove our main theorem.

\section{Proof of Main Theorem}

\begin{figure}
\centering
\begin{tikzpicture}[/pgf/declare function={f=8/x;}]
\begin{axis}[
        xmin=0,xmax=9,ymin=0,ymax=4,
        xtick={0,1,2.75,4.5,6.25,8},
    xticklabels={$0$,$e^{-4}$,$e^{-3}$,$e^{-2}$, $e^{-1}$, 1},
    ytick={0,7.1,7.2,7.3,7.5,8},
    yticklabels={$$,$$,$$,$$, $$, 1},
    domain=0:10,
    samples=100,
    axis lines=middle
]
\addplot [thick, red] {f}
    node [pos=.9, above right] {$f_n(x)$};

\addplot [
    black!80,fill=green,opacity=.3,
    right segments=4,
    right=1:8,
] {f};
\end{axis}
\end{tikzpicture}
\caption{Right endpoint rule approximation for $f_n(x)$.} \label{fig1}
\end{figure}
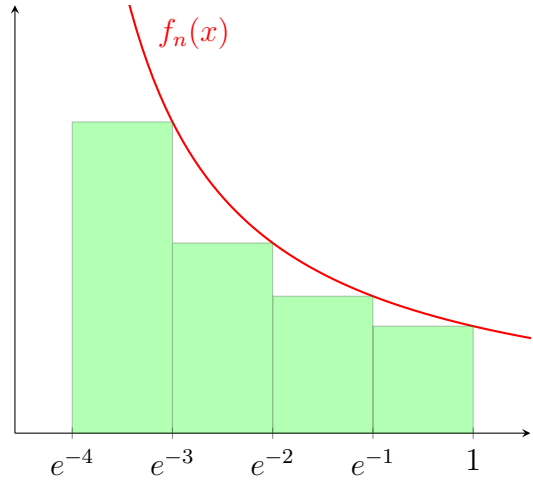

\begin{atmProof}
Recall that $p$ to denotes a prime number. We will be working with the following series.
\begin{equation*}
f_p(x) = \sum_{k = 0}^{\infty} x^{\log(p^k)} = \sum_{k = 0}^{\infty} x^{k \log(p)} = \frac{1}{1 - x^{\log(p)}}
\end{equation*}
which converges for $0 \le x < 1.$ The Fundamental Theorem of Arithmetic tells us that every integer less than or equal to $n$ can only use primes less than or equal to $n$ in its prime factorization. This gives us the following inequality:

\begin{equation}\label{eq1}
\prod_{p \le n} f_p(x) \ge 1 + x^{\log(2)} + x^{\log(3)} + \cdots x^{\log(n)} = \sum_{k = 1}^{n} x^{\log(k)}.
\end{equation}
Now let $$F_n(x) = \sum_{k = 1}^{n} x^{\log(k)}.$$
Inequality \eqref{eq1} gives us 
$$\sum_{p \le n} -\log(1 - x^{\log(p)}) \ge \log(F_n(x))$$
and therefore,
\begin{equation} \label{eq2}
    \int_0^1 \sum_{p \le n} -\frac{ \log(1 - x^{\log(p)})}x \ dx \ge \int_0^1 \frac{\log(F_n(x))}x \ dx
\end{equation}
Exchanging the sum and the integral (since our sum is finite) in the left hand side of Inequality \eqref{eq2} gives us
\begin{equation*}
\sum_{p \le n} \int_0^1 - \frac{ \log(1 - x^{\log(p)})}x \ dx  \ge \int_0^1 \frac{\log(F_n(x))}x \ dx.
\end{equation*}
We now make the substitution $y = x^{\log(p)}$. 
\begin{equation*}
\sum_{p \le n}\frac{1}{\log(p)} \int_0^1 - \frac{ \log(1 - y)}y \ dy  \ge \int_0^1 \frac{\log(F_n(x))}x \ dx.
\end{equation*}
Recall the well known identity
$$\int_0^1 - \frac{\log(1-y)}y \ dy = \frac{\pi^2}6$$
which follows from the Taylor expansion of $\frac{\log(1-y)}y$ and using the fact that $\sum_{n \ge 1} \frac{1}{n^2} = \frac{\pi^2}6.$ Thus,
\begin{equation}\label{eq3}
\frac{\pi^2}6 \sum_{p \le n} \frac{1}{\log(p)} \ge \int_0^1 \frac{\log(F_n(x))}x \ dx.
\end{equation}
Since $\log(x) = o(x)$, we can check that $$f_n(x) = \frac{\log(F_n(x))}x$$ is monotonically decreasing over the interval $[0, \infty)$ for $n \ge 2$. Therefore, we can bound the integral in the right hand side of Inequality \eqref{eq3} using the right endpoint rule. (See Figure \ref{fig1}.) This gives us
\begin{align*}
\int_0^1 f_n(x) \ dx &> \sum_{k = 0}^{\infty} (e^{-k} - e^{-k-1}) f_n(e^{-k})  \\
&= \sum_{k=0}^{\infty} (1-e^{-1}) \log \left( \sum_{j = 1}^n j^{-k} \right) \\
&= (1-e^{-1})\left( \log(n) +  \log(H_n) \right) + (1-e^{-1})\sum_{k=2}^{\infty} \log\left( \sum_{j = 1}^n j^{-k} \right)
\end{align*}
where 
\begin{equation} \label{eq5}
H_n =  \sum_{j=1}^n \frac{1}j \sim \log(n).
\end{equation}
We will now bound the term $\log\left( \sum_{j = 1}^n j^{-k} \right)$ . Let $\zeta(k)$ denote the Riemann zeta function defined for $s \in \mathbb C$ by the sum $$\zeta(s) := \sum_{n=1}^{\infty} \frac{1}{n^s},\ \ \operatorname{Re}(s)>1.$$ We have 
\begin{align}
\sum_{k=2}^{\infty} \log\left( \sum_{j = 1}^n j^{-k} \right) &\le \sum_{k=2}^{\infty} \log(\zeta(k)) 
= \sum_{k = 2}^{\infty} \log(1 + (\zeta(k)-1)) \nonumber \\
&< \sum_{k=2}^{\infty} (\zeta(k)-1)
= \sum_{k=2}^{\infty} \sum_{n = 2}^{\infty} \frac{1}{n^k}\nonumber \\
&= \sum_{n=2}^{\infty} \sum_{k=2}^{\infty} \frac{1}{n^k} = \sum_{n=2}^{\infty} \left(\frac{1}{n-1} - \frac{1}n \right) = 1. \label{eq6}
\end{align}  
Using the results of Equations \eqref{eq5} and \eqref{eq6} in our computation of the integral of $f_n(x)$ tells us
\begin{equation*}
\int_0^1 f_n(x) \ dx > (1-e^{-1}) (\log(n) + \log\log(n)) + A
\end{equation*}
where $A$ is a positive constant less than $1-e^{-1}.$ 
Inequality \eqref{eq3} gives us 

\begin{align*}
    \sum_{p \le n} \frac{1}{\log(p)} \ge \frac{6}{\pi^2} \int_0^1 f_n(x) \ dx > \frac{6(1-e^{-1})}{\pi^2}\log(n) > \frac{1}3 \log(n)
\end{align*}
as desired. 
\end{atmProof}


\bibliographystyle{plain}
\bibliography{main}

\begin{thebibliography}{1}

\bibitem{thebook}
Martin Aigner, Ziegler~G\"unter M., and Karl~H. Hofmann.
\newblock {\em Proofs from THE BOOK}.
\newblock Springer, 2018.

\bibitem{euler}
William Dunham.
\newblock {\em Euler: the master of us all}.
\newblock Mathematical Association of America, 1999.

\bibitem{rosen}
Kenneth~F. Ireland and Michael~I. Rosen.
\newblock {\em A classical introduction to modern number theory}.
\newblock Springer-Verlag, 2010.

\bibitem{PNT}
D.~Zagier.
\newblock Newman's short proof of the prime number theorem.
\newblock {\em The American Mathematical Monthly}, 104(8):705--708, 1997.

\end{thebibliography}

\end{document}